\documentclass[a4paper]{article}
\usepackage{amsmath,amssymb,amsthm,amscd}

\newtheorem{Thm}{Theorem}[subsection]
\newtheorem{Lem}[Thm]{Lemma}
\newtheorem{Cor}[Thm]{Corollary}
\newtheorem{Prop}[Thm]{Proposition}
\newtheorem{Def}[Thm]{Definition}
\newtheorem{Rem}[Thm]{Remark}

\makeatletter
\@addtoreset{equation}{section}
\makeatother

\newcommand{\R}{\mathbb{R}}

\newcommand{\N}{\mathbb{N}}

\newcommand{\Z}{\mathbb{Z}}

\newcommand{\Sch}{\mathcal{S}}
\newcommand{\E}{\mathcal{E}}
\newcommand{\Li}{\mathcal{L}}

\newcommand{\Fu}{\mathcal{F}}

\newcommand{\Rb}{\overline{\R}}

\newcommand{\dbar}{{d\hspace{-0,05cm}\bar{}\hspace{0,05cm}}}
\newcommand{\op}{\textup{Op}}

\newcommand{\clw}{\textup{cl}}

\newcommand{\cone}{\textup{cone}}

\newcommand{\la}{\langle}
\newcommand{\ra}{\rangle}
\newcommand{\lr}{\textup{(}}
\newcommand{\rr}{\textup{)}}
\newcommand{\diff}{\textup{Diff}}

\newcommand{\td}{\tilde}

\date{}
\title{Edge-degenerate families of $\Psi$DO's on an infinite cylinder}
\author{Jamil Abed \and Bert-Wolfgang Schulze}
\begin{document}
\maketitle
\begin{abstract} 
We establish a parameter-dependent pseudo-differential calculus on an infinite cylinder, regarded as a manifold with conical exits to infinity.
The parameters are involved in edge-degenerate form, and we formulate the operators in terms of operator-valued amplitude functions.\\
\linebreak
\end{abstract}
\tableofcontents
\section*{Introduction}
\addcontentsline{toc}{section}{Introduction}
\markboth{INTRODUCTION}{INTRODUCTION}
The analysis of (pseudo-) differential operators on a manifold (stratified space) with higher polyhedral singularities employs to a large extent parameter-dependent families of operators on a (in general singular) base $X$ of a cone, where the parameters $(\rho,\eta)$ have the meaning of covariables in cone axis and edge direction, respectively. These covariables appear in edge-degenerate form, i.e., in the combination $(r\rho,r\eta)$ where $r\in\R_+$ is the axial variable of the cone $X^\Delta=(\Rb_+\times X)/(\{0\}\times X)$ with base $X$. For operators on a wedge $X^\Delta\times\Omega\ni (\cdot,y)$ it is essential to understand the structure of what we call edge symbols, together with associated weighted distributions on the cone. It is natural to split up the investigation into a part for $r\to 0$, i.e., close to the tip of the cone and a part for $r\to\infty$, the conical exit to infinity. For higher corner theories it is desirable to do that in an axiomatic manner, i.e., to point out those structures which make the calculus iterative. For the case $r\to 0$ the authors developed in \cite{Abed1} such an axiomatic approach. What concerns $r\to\infty$ it seems to be advisable first to concentrate on the case when the base $X$ is smooth and compact. It turns out that the edge-degeneration of symbols in a calculus on $\R$ up to infinity causes a highly non-standard behaviour with respect to symbolic rules for operator-valued symbols (to be invented in the right manner). This has to be analysed first, where the approach should rely on the principles and key properties that are essential for the iteration. The goal of our paper is just to develop some crucial steps in that sense.\\
\indent To be more precise, we show (here for a smooth compact manifold as the base of a cone) how very simple and general phenomena on the norm-growth of parameter-dependent pseudo-differential operators in the sense of Theorem \ref{131.red2707.th} are sufficient to induce the essential properties of a calculus on the manifold $X^\asymp\cong \R\times X$ with conical exits $r\to\pm \infty$. In other words, knowing a suitable variant of Theorem \ref{131.red2707.th} for a singular (compact) manifold, we can expect essentially the same things on the respective singular cylinder. Details in that case go beyond the scope of the present paper; let us only note that the article \cite{Abed1} just contains an analogue of Theorem \ref{131.red2707.th} for a base manifold with edge. We introduce here a pseudo-differential calculus in spaces $H^{s,g}_\cone(X^\asymp)$ in a self-contained manner, including those spaces themselves. In the smooth case there is, of course, also a completely independent approach, usually organised without parameter $\eta$, well-known under the key-words operators on manifolds with conical exit to infinity, here realised on such a manifold $X^\asymp$ modelled on an infinite cylinder. Concerning the generalities we refer to Shubin \cite{Shub1}, Parenti \cite{Pare1}, Cordes \cite{Cord1}, or to the corresponding sections in \cite{Schu20}. There are also several papers for singular $X$, cf., for instance, \cite{Schu25}, or \cite{Calv1}, \cite{Calv3}, but those are based on more direct information from corner-degenerate symbols. We think that the present idea admits to manage the iterative process for higher singularities in a more transparent way.

\section{A new class of operator-valued symbols}
\label{cylsec1}
\subsection{Edge-degenerate families on a smooth compact manifold}
Edge-degenerate families of pseudo-differential operators occur in connection with the edge symbols of operators of the form
\begin{equation} \label{cylA}
A=r^{-\mu}\sum_{j+|\alpha|\leq\mu}a_{j\alpha}(r,y)\left(-r\partial_r\right)^j(rD_y)^\alpha
\end{equation}
with coefficients $a_{j\alpha}\in C^\infty\big(\Rb_+\times\Omega,\diff^{\mu-(j+|\alpha|)}(X)\big)$ for an open set $\Omega\subseteq\R_y^q$; here $X$ is a smooth compact manifold, and $\diff^\nu(X)$ is the space of all differential operators of order $\nu$ on $X$ with smooth coefficients. The analysis of such edge-degenerate operators is crucial for understanding the solvability of elliptic equations on spaces with polyhedral singularities, (cf. \cite{Schu27}, \cite{Calv2}, or \cite{Mani2}). Apart from the standard homogeneous principal symbol of \eqref{cylA} which is a $C^\infty$ function on $T^*(\R_+\times X\times\Omega)\setminus 0$, we have the so-called principal edge symbol
\begin{equation} \label{cyleds}
\sigma_\wedge(A):=r^{-\mu}\sum_{j+|\alpha|\leq\mu}a_{j\alpha}(0,y)\left(-r\partial_r\right)^j(r\eta) ^\alpha
\end{equation}
parametrised by $(y,\eta)\in T^*\Omega\setminus 0$ and with values in Fuchs type differential operators on the open infinite stretched cone $X^\wedge:=\R_+\times X$ with base $X$. For the construction of parametrices of $A$ (in the elliptic case) we need to understand, in particular, the nature of parameter-dependent parametrices of operator families 
\begin{equation} \label{cyldeg}
r^{-\mu}\sum_{j+|\alpha|\leq\mu}a_{j\alpha}(0,y)(-ir\rho)^j(r\eta)^\alpha
\end{equation}
on $\R_+\times X$ for $r\to\infty$. We often set $\td\rho=r\rho,\td\eta=r\eta$. If $A$ is edge-degenerate elliptic (cf. \cite{Schu20},\cite{Haru13}) it turns out that $\sum_{j+|\alpha|\leq\mu}a_{j\alpha} (0,y)(-i\td\rho)^j(\td\eta)^\alpha$ is parameter-dependent elliptic on $X$ with parameters $(\td\rho, \td\eta)\in\R^{1+q}$, for every fixed $y\in\Omega$. Let $L_\clw^\mu(X;\R^l)$ denote the set of all parameter-dependent classical pseudo-differential operators of order $\mu\in\R$ on the manifold $X$, with parameters $\lambda\in\R^l,l\in\N$. That means, the amplitude functions $a(x,\xi,\lambda)$ in local coordinates $x\in\R^n$ on $X$ are classical symbols of order $\mu$ in $(\xi,\lambda)$. The space $L^{-\infty}(X;\R^l)$ of parameter-dependent smoothing operators is defined via kernels in $\Sch\big(\R^l,C^\infty(X\times X)\big)$ (a fixed Riemannian metric on $X$ admits to identify $C^\infty(X\times X)$ with corresponding integral operators).\\
For future references we state and prove a standard property of parameter-dependent operators.
\begin{Thm} \label{131.red2707.th}
Let $M$ be a closed compact $C^\infty$ manifold and $A(\lambda)\in L_\clw^\mu(M;\R^l)$ a parameter-dependent family of order $\mu$, and let $\nu\geq\mu$. Then there is a constant $c=c(s,\mu,\nu)>0$ such that 
\begin{equation} \label{cyleqlam}
\|A(\lambda)\|_{\Li(H^s(M),H^{s-\nu}(M))}\leq c\la\lambda\ra^{\max\{\mu,\mu-\nu\}}.
\end{equation}
In particular, for $\mu\leq 0$, $\nu=0$ we have
\begin{equation} \label{cyleq2}
\|A\|_{\Li(H^s(M),H^{s}(M))}\leq c\la\lambda\ra^\mu.
\end{equation}
Moreover, for every $s',s''\in\R$ and every $N\in\N$ there exists a $\mu(N)\in\R$ such that for every $\mu\leq\mu(N)$, $k:=\mu(N)-\mu$, and $A(\lambda)\in L_\clw^\mu(M;\R^l)$ we have
\begin{equation} \label{cyleq3}
\|A\|_{\Li(H^{s'}(M),H^{s''}(M))}\leq c\la\lambda\ra^{-N-k}.
\end{equation}
for all $\lambda\in\R^l$, and a constant $c=c(s',s'',\mu,N,k)>0$.
\end{Thm}
\proof
In this proof we write $\|\cdot\|_{s',s''}=\|\cdot\|_{\Li(H^{s'}(M),H^{s''}(M))}$. The estimates \eqref{cyleqlam} and \eqref{cyleq2} are standard. Concerning \eqref{cyleq3} we first observe that we have to choose $\mu$ so small that $A(\lambda):H^{s'}(M)\to H^{s''}(M)$ is continuous. This is the case when $s''\leq s'-\mu$, i.e., $\mu\leq s'-s''$. Let $R^{s''-s'}(\lambda)\in L_\clw^{s''-s'} (M,\R^l)$ be an order reducing family with the inverse $R^{s'-s''}(\lambda)\in L_\clw^{s'-s''}(M, \R^l)$. Then we have 
\[ R^{s''-s'}(\lambda):H^{s''}(M)\to H^{s'}(M), \]
i.e., $R^{s''-s'}(\lambda)A(\lambda):H^{s'}(M)\to H^{s'}(M)$. The estimate \eqref{cyleq2} gives us
\[ \|R^{s''-s'}(\lambda)A(\lambda)\|_{s',s'}\leq c\la\lambda\ra^{\mu+(s''-s')} \]
for $\mu\leq s'-s''$. Moreover, \eqref{cyleqlam} yields $\|R^{s'-s''}(\lambda)\|_{s',s''}\leq c\la\lambda\ra^{s'-s''}$. Thus 
\begin{multline*} \|A(\lambda)\|_{s',s''}=\|R^{s'-s''}(\lambda)R^{s''-s'}(\lambda)A(\lambda)\|_{s',s''} \\ \leq \|R^{s'-s''}(\lambda)\|_{s',s''}\|R^{s''-s'}(\lambda)A(\lambda)\|_{s',s'}\leq c\la\lambda\ra ^{(s'-s'')+\mu+(s''-s')}=c\la\lambda\ra^\mu.
\end{multline*}
In other words, when we choose $\mu(N)$ in such a way that $\mu\leq {s'-s''}$, and $\mu(N)\leq -N$, then \eqref{cyleq3} is satisfied. In addition, if we take $\mu=\mu(N)-k$ for some $k\geq 0$ then \eqref{cyleq3} follows in general.
\qed
\begin{Cor} \label{131.spec2707.co}
Let $A(\lambda)\in L_\clw^\mu(M;\R^l)$, and assume that the estimate
\[ \|A(\lambda)\|_{s',s''}\leq c\la\lambda\ra^{-N} \]
is true for given $s',s''\in\N$ and some $N$. Then we have
\[ \|D_\lambda^\alpha A(\lambda)\|_{s',s''}\leq c\la\lambda\ra^{-N-\alpha} \]
for every $\alpha\in\N^l$.
\end{Cor}

Since we are interested in families for $r\to\infty$ it will be convenient to ignore the specific edge-degenerate behaviour for $r\to 0$ and to consider the cylinder $\R\times X$ rather than $\R_+\times X$. Far from $r=\pm\infty$ our calculus will be as usual; therefore, for convenience, we fix a strictly positive function $r\to[r]$ in $C^\infty(\R)$ such that $[r]=|r|$ for $|r|>R$ for some $R>0$. The operator-valued amplitude functions in our calculus on $\R\times X$ are operator families of the form 
\[a(r,\rho,\eta)=\td a(r,[r]\rho,[r]\eta) \]
where $\td a(r,\td\rho,\td\eta)\in C^\infty\big(\R,L_\clw^\mu(X;\R^{1+q}_{\td\rho,\td\eta})\big)$. In addition it will be important to specify the dependence of the latter function for large $|r|$. In other words, the crucial definition is as follows.
\begin{Def}
\begin{enumerate}
\item[\textup{(i)}] Let $E$ be a Fr\'echet space with the \lr countable\rr\ system of semi-norms $(\pi_j) _{j\in\N}$; then $S^\nu(\R,E)$, $\nu\in\R$, is defined to be the set of all $a(r)\in C^\infty(\R,E)$ such that
\[ \pi_j\left(D_r^ka(r)\right)\leq c[r]^{\nu-k} \]
for all $r\in\R,k\in\N$, with constants $c=c(k,j)>0$,
\item[\textup{(ii)}]$\boldsymbol{S}^{\mu,\nu}$ for $\mu,\nu\in\R$ denotes be the set of all operator families
\[ a(r,\rho,\eta)=\td a(r,[r]\rho,[r]\eta) \]
for $\td a(r,\td\rho,\td\eta)\in S^\nu\big(\R,L_\clw^\mu(X;\R^{1+q}_{\td\rho,\td\eta})\big)$ \lr referring to the natural nuclear Fr\'echet topology of the space $L_\clw^\mu(X;\R^{1+q} _{\td\rho,\td\eta})$\rr.
\end{enumerate}
\end{Def}
We first establish some properties of $\boldsymbol{S}^{\mu,\nu}$ that play a role in our calculus.
\begin{Prop}
\label{152.(cal)2209.pr}
\begin{enumerate} 
\item[\textup{(i)}] $\varphi(r) \in S^\sigma (\R),a(r, \rho, \eta) \in
\boldsymbol{S}^{\mu,\nu}$ implies $\varphi(r)a(r,\rho,\eta)$ $\in\boldsymbol{S}^{\mu,\sigma+\nu}$.
\item[\textup{(ii)}] For every $k,l \in \N$ we have 
   \[ a \in \boldsymbol{S}^{\mu,\nu} \Rightarrow \partial_r^l a \in
       \boldsymbol{S}^{\mu,\nu-l}, 
         \partial_\rho^k a \in \boldsymbol{S}^{\mu-k,\nu+k},
	 \partial_\eta^k a \in \boldsymbol{S}^{\mu-k,\nu+k}.   \]
\item[\textup{(iii)}] $a(r,\rho,\eta)\in\boldsymbol{S}^{\mu,\nu}$, $b(r,\rho,\eta)\in\boldsymbol{S} ^{\td\mu,\td\nu}$ implies $a(r,\rho,\eta)b(r,\rho,\eta)\in\boldsymbol{S}^{\mu+\td\mu,\nu+\td\nu}$.
\end{enumerate}	 
\end{Prop}

\proof
(i) is evident. (ii) For simplicity we assume 
    $q = 1$ and compute
\[  \partial_r\td a(r,[r]\rho,[r]\eta) = 
           \big((\partial_r + [r]'\rho
	   \partial_{\td\rho} + [r]'\eta
	   \partial_{\td\eta})\td a \big)
	   (r,[r]\rho,[r]\eta)    \]
    where $[r]' := \partial_r[r]$. Since
    $\td\rho\td a(r, \td\rho,\td\eta)$,
        $\td\eta \td a(r, \td\rho, 
	\td\eta) \in S^\nu(\R,
        L_\clw^{\mu+1}(X;\R_{\td\rho,
	\td\eta}^{1+q}))$, and
 	$\partial_{\td\rho} \td a,
	\partial_{\td\eta} \td a \in S^\nu(\R,
	L_\clw^{\mu-1}(X; \R^{1+q}))$, we obtain  
    \[ \partial_r \td a(r,[r]\rho, [r]\eta)
         = \big((\partial_r+([r]'/[r])[r]\rho
	 \partial_{\td\rho}+ ([r]'/[r])[r] \eta
	 \partial_{\td\eta})\td a\big)
	    (r,[r]\rho, [r]\eta) \in \boldsymbol{S}^{\mu,\nu-1}.  \]   
It follows that $\partial_r^l a \in \boldsymbol{S}^{\mu,\nu-l}$ for all $l \in \N$.
Moreover, we have
\[  \partial_\rho \td a(r,[r]\rho, [r]\eta) = 
       [r](\partial_{\td\rho} \td a) (r,[r]\rho,
       [r] \eta)   \]
which gives us $\partial_\rho a \in \boldsymbol{S}^{\mu-1,
\nu+1}$, 
and,
by iteration, $\partial_\rho^k a \in \boldsymbol{S}^{\mu-k,\nu+k}$. In
a similar manner we can argue for the $\eta$-derivatives.\\ 
(iii) By definition we have
\[ a(r,\rho,\eta)=\td a(r,[r]\rho,[r]\eta), \ b(r,\rho,\eta)=\td b(r,[r]\rho,[r]\eta) \]
for $\td a(r,\td\rho,\td\eta)\in S^\nu\big(\R,L_\clw^\mu(X,\R^{1+q}_{\td\rho,\td\eta})\big)$, $\td b(r,\td\rho,\td\eta)\in S^{\td\nu}\big(\R,L_\clw^{\td\mu}(X,\R^{1+q}_{\td\rho,\td\eta})\big)$. Then the assertion is a consequence of the relation
\[ (\td a\td b)(r,\td\rho,\td\eta)\in S^{\nu+\td\nu}\big(\R,L_\clw^{\mu+\td\mu}(X,\R^{1+q}_ {\td\rho,\td\eta})\big). \]
\qed 	

\begin{Cor} \label{cylcormul}
For $a(r,\rho,\eta)\in\boldsymbol{S}^{\mu,\nu}$, $b(r,\rho,\eta)\in\boldsymbol{S}^{\td\mu,\td\nu}$ for every $k\in\N$ we have 
\[ \partial_\rho^k a(r,\rho,\eta) D_r^k b(r,\rho,\eta)\in\boldsymbol{S}^{\mu+\td\mu-k,\nu+\td\nu} \]
\end{Cor}
\proof
In fact, we have $\partial_\rho^k a(r,\rho,\eta)\in\boldsymbol{S}^{\mu-k,\nu+k}$, $\partial_r^k b(r,\rho,\eta) \in\boldsymbol{S}^{\td\mu,\td\nu-k}$.
\qed
\begin{Prop}
\label{152.2507.pr}
Let $\td a_j (r,\td\rho,\td\eta) \in S^\nu(\R,
L_\clw^{\mu-j}(X; \R^{1+q}))$, $j \in \N$, be an  arbitrary
sequence, $\mu, \nu \in \R$ fixed. Then there is an
$\td a(r,\td\rho,\td\eta) \in S^\nu(\R,
L_\clw^\mu(X; \R^{1+q}))$ such that
\[  a- \sum_{j=0}^{N} a_j \in S^\nu(\R, L_\clw^{\mu-(N+1)} (X;
      \R^{1+q}))     \]
for every $N \in \N$, and a is unique $\bmod\, S^\nu \big(\R, L_\clw^{-
\infty}(X; \R^{1+q})\big)$.
\end{Prop}      

\proof
The proof is similar to the standard one on asymptotic summation of
symbols. We can find
an asymptotic sum as a convergent series $\td a
(r,\td\rho,\td\eta) = \sum_{j=0}^{\infty}
       \chi \left( (\td\rho,\td\eta)/c_j\right)
       \td a_j(r,\td\rho,\td\eta)$  
for some excision function $\chi$ in $\R^{1+q}$, with a sequence
$c_j > 0$, $c_j \to \infty$ as $j \to \infty$ so fast, that
$\sum_{j=N+1}^{\infty} \chi \left(
(\td\rho,\td\eta)/c_j\right)\td a(r,\td\rho,
\td\eta)$ converges in $S^\nu(\R, L_\clw^{\mu-(N+1)})$ for
every $N$.       
\qed

\subsection{Continuity in Schwartz spaces}
\begin{Thm}
\label{152.2309.th}
Let $p(r,\rho,\eta) = \td p(r,[r]\rho,[r]\eta)$,  
$\td p(r,\td\rho,\td\eta) \in S^\nu\big(\R,
L_\clw^\mu(X;$ $\R_{\td\rho,\td\eta}^{1+q})\big)$,
i.e., $p(r,\rho,\eta) \in \boldsymbol{S}^{\mu,\nu}$. Then
$\op_r(p)(\eta)$ induces a family of continuous operators
\[  \op_r(p)(\eta) : {\Sch}\big(\R, C^\infty(X)\big) \to \Sch\big(\R, C^\infty(X)\big)   \]
for every $\eta\not=0$.
\end{Thm}

\proof
We have
\[  \op_r(p)(\eta) u(r) = \int e^{ir \rho}
      p(r,\rho,\eta) \hat{u}(\rho) \dbar \rho,   \]
first for $u \in C_0^\infty\big(\R, C^\infty(X)\big)$. 
In the space ${\Sch}\big(\R, C^\infty(X)\big)$ we have
the semi-norm system
\[   \pi_{m,s}(u) = \max_{\alpha+\beta \leq m} \sup_{r \in
       \R} \| [r]^\alpha \partial_r^\beta u(r) \|_{H^s (X)}
       \] 
for $m \in \N$, $s \in \Z$, which defines the Fr\'echet
topology of ${\Sch} \big(\R, C^\infty(X)\big)$. 

If necessary we indicate the variable $r$, i.e., write
$\pi_{m,s;r}$ rather than $\pi_{m,s}$. The Fourier transform
$\Fu_{r \to \rho}$ induces an isomorphism
\[ \Fu: {\Sch}\big(\R_r, H^s(X)\big) \to  {\Sch}\big(\R_\rho, H^s(X)\big)    \]
for every $s$. For every $m \in \N$ there exists a $C > 0$
such that 
\begin{equation}
\label{152.(sem)2309.eq}
\pi_{m,s;\rho} (\Fu u) \leq C \pi_{m+2,s;r} (u)
\end{equation}
for all $u \in {\Sch}(\R, H^s(X))$ (see
\cite[Chapter 1]{Kuma1} for scalar functions; the case of
functions with values in a Hilbert space is completely
analogous).      
We have to show that for every $\td m \in \N$, $\td s
\in \Z$ there exist $m \in \N$, $s \in \Z$, such that
\begin{equation}
\label{152.(est)2309.eq}
\pi_{\td m,\td s} \big((\op(p) u)(r)\big) \leq c \pi_{m,s}(u)
\end{equation}
for all $u \in {\Sch}\big(\R,C^\infty(X)\big)$, for some $c =
c(\td m,\td s) > 0$. According to Proposition
\ref{152.2209.co} below we write the operator $\op(p)(\eta)$ in the
form
\begin{equation}
\label{152.(sum)2309.eq}
\op_r(p)(\eta) \circ \la r \ra^{-M} \circ \la r
    \ra^M = \la r \ra^{-M} \op_r(b_{MN})(\eta) 
    \circ \la r \ra^M + \op_r(d_{MN})(\eta) \circ \la r
    \ra^M
\end{equation}
for a symbol $b_{MN}(r, \rho,\eta) \in \boldsymbol{S}^{\mu,\nu}$ and
a remainder $d_{MN}(r,\rho,\eta)$ satisfying  
estimates analogously as \eqref{152.(NE)1608.eq}.    

We have
\begin{multline}
\label{152.(ST)2309.eq}
\| \op_r(p)(\eta)  u(r) \|_{H^{\td s}(X)} 
  = \| \int e^{ir \rho} p(r,\rho,\eta) \hat{u}(\rho)
       \dbar \rho \|_{H^{\td s}(X)}    \\ 
  \leq \| \int e^{ir \rho} \la r \ra^{- M}
               b_{MN}(r,\rho,\eta)\big( \la r \ra^M
	       u\big)^\land (\rho) \dbar \rho
	       \|_{H^{\td s}(X)} \\
	       + \| \op_r(d_{MN})(\eta) (\la r \ra^M u(r))
	       \|_{H^{\td s}(X)}.
\end{multline}
For the first term on the right of \eqref{152.(ST)2309.eq} we
obtain for $s := \td s + \mu$ and arbitrary $\widetilde M \in\N$
\begin{multline*}
\| \int e^{ir \rho} \la \rho
   \ra^{- \widetilde{M}} \la r \ra^{- M} b_{MN}
   (r,\rho, \eta) \la \rho \ra^{\widetilde{M}}
   \big(\la r \ra^M u\big)^\land (\rho) \dbar \rho
   \|_{H^{\td s}(X)}   \\	            
\leq  \int \| \la \rho \ra^{- \widetilde{M}}
    \la r \ra^{- M} b_{MN}(r,\rho,\eta) \la
    \rho \ra^{\widetilde{M}} \big(\la r \ra^M
    u \big)^\land(\rho) \|_{H^{\td s}(X)} \dbar \rho
         \\  
\leq c \sup_{(r, \rho)\in \R^2} \la \rho
         \ra^{- \widetilde{M}} \la r \ra^{- M} 
	 \| b_{MN}(r,\rho,\eta) \|_{ {\Li}(H^s(X),H^{\td s}(X))} \\
     \int \la \rho \ra^{\widetilde{M}} \|
       \big(\la r \ra^M u\big)^\land (\rho) \|_{H^s(X)} \dbar\rho.
\end{multline*}

Moreover, we have        	  
\begin{align*}
\int \la \rho \ra^{\widetilde{M}} \| \big( \la r
     \ra^M & u\big)^\land (\rho) \|_{H^s(X)} \dbar \rho   
\leq  \sup_{\rho \in \R} \la \rho
         \ra^{\widetilde{M}+2} \| \big( \la r 
	 \ra^M u\big)^\land (\rho) \|_{H^s(X)}
	 \int \la \rho \ra^{- 2} \dbar \rho \\
\leq & \ c \pi_{\widetilde{M}+2,s;\rho} \big(\big( \la r
         \ra^M u\big)^\land (\rho)\big)   
\leq  \pi_{\widetilde{M}+ 4,s; r} (\la r \ra^M u))
\leq c \pi_{M+\widetilde{M}+ 4,s;r} (u)
\end{align*}
Here we employed the estimate \eqref{152.(sem)2309.eq}.	 
Thus \eqref{152.(ST)2309.eq} yields
\begin{multline}
\label{152.(es1)2308.eq}
\pi_{0,\td s}(\op(p)(\eta) u) \leq c \sup_{(r,\rho)\in \R^2} \la \rho \ra^{- \widetilde{M}} \la
     r \ra^{- M} \| b_{MN} (r, \rho, \eta) \|_{ {\Li}(H^s(X), H^{\td s}(X))} \\
  \pi_{M+\widetilde{M}+4,s;r}(u) + \pi_{0,\td s} (\op_r(d_{MN})(\eta) ( \la r \ra^M
     u)).
\end{multline}
The factor $c \sup_{(r, \rho)\in \R^2} \la \rho\ra^{- \widetilde{M}} \la r \ra^{- M}\|  b_{MN}(r,\rho,\eta) \|_{ {\Li}(H^s(X),H^{\td s}(X))}$ in front of $\pi_{M+\widetilde{M}+4}(u)$ is finite when we choose $M$ so large that $\nu - M \leq 0$ and
$\widetilde{M}$ so large that
\[  \sup_{\rho \in \R} \la \rho \ra^{-
      \widetilde{M}} \| b_{MN}(r, \rho,\eta) \|_ { {\Li}(H^s(X), H^{\td s}(X))} < \infty. \]
Next we consider the second term on the right hand side of
\eqref{152.(es1)2308.eq}. We have
\begin{multline*}
\| \op_r(d_{MN})(\eta) \la r \ra^M u(r)\|_{H^{\td s}(X)} \\
=\left\| \int e^{ir \rho} \la\rho \ra^{-M} d_{MN}(r, \rho,\eta) \la \rho\ra^M (\la r \ra^M u)^\land(\rho)\dbar \rho\right\|_{H^{\td s}(X)}    \\
\leq \int \| \la \rho \ra^{- M} d_{MN} (r,\rho,\eta) \la \rho \ra^M \big(\la r \ra^M
         u\big)^\land (\rho) \|_{H^{\td s}(X)} \dbar \rho \\
\leq \int \sup_{(r,\rho)\in \R^2} \| \la \rho\ra^{- M} d_{MN}(r, \rho, \eta) \|_{ {\Li}(H^{\td s}(X), H^s(X))}\|\la \rho \ra^M \big(\la r \ra^M
	  u\big)^\wedge(\rho) \|_{H^s(X)} \dbar \rho.
\end{multline*}

From the analogue of  the estimate \eqref{152.(NE)1608.eq}
for $d_{MN}(r,\rho,\eta)$ we see that for $N$ sufficiently
large it follows that the right hand side of the latter
expression can be estimated by
\begin{align*}
c \int \| \la \rho \ra^M \big( \la r \ra^M
& u\big)^\land(\rho)  \|_{H^s(X)} \dbar \rho \\
& \leq \sup_{\rho\in\R} \la \rho \ra^{M+2} \| \big(\la r \ra^M
    u\big)^\land  (\rho) \|_{H^s(X)} \int \la \rho
    \ra^{-2} \dbar \rho. \\	
& \leq c \pi_{2M+2,s;\rho} (\hat{u}(\rho)) \leq c
     \pi_{2M+4,s;r} (u).
\end{align*}     

In other words we proved 
\begin{equation}
\label{152.(fin)2409.eq}
\pi_{0,\td s} (\op(p)(\eta)) \leq c \{\pi_{M+\widetilde{M}+4,s} (u)+\pi_{2M+4,s}(u)\} \leq c \pi_{L,s} (u)
\end{equation}
for $s = \td s+\mu$, $L := \max \{ M+\widetilde{M}+4,2M+4
\}$.
Now we write
\begin{align*}
\partial_r \op(p)(\eta) u(r) & = \int e^{ir \rho} \{ i
\partial_\rho + \partial_r \} p(r,\rho,\eta) \hat{u}(\rho) \dbar\rho,   \\
r \op(p)(\eta)u(r) & = \int e^{ir \rho} (i \partial_\rho
p(r,\rho,\eta)) \hat{u}(\rho) \dbar \rho+\int e^{ir\rho}ip(r,\rho,\eta)\partial_\rho\hat u(\rho)\dbar\rho.		        	  	 	
\end{align*}				
From Proposition \ref{152.(cal)2209.pr} we have
\[  \{ i \partial_\rho + \partial_r \} p(r,\rho,\eta) \in
      \boldsymbol{S}^{\mu-1,\nu+1} + \boldsymbol{S}^{\mu,\nu-1} \subseteq
      \boldsymbol{S}^{\mu,\nu+1}, \
      i \partial_\rho p(r,\rho,\eta) \in
      \boldsymbol{S}^{\mu-1,\nu+1}.     \]
Since the estimate \eqref{152.(fin)2409.eq} is  true for elements
in the respective symbol classes of arbitrary order, it follows
altogether the estimate \eqref{152.(est)2309.eq} for every
$\td m\in\N ,\td s\in\Z$ and suitable $m,s$.      
\qed

\subsection{Leibniz products and remainder estimates}
\label{cylsecleib}
Let $\td a(r,\td\rho, \td\eta) \in S^\nu(\R,
     L_\clw^\mu),\td b(r, \td\rho,\td\eta) \in
     S^{\td\nu}(\R, L_\clw^{\td\mu})$
where $L_\clw^\mu = L_\clw^\mu(X; \R_{\td\rho,
\td\eta}^{1+q})$.
The operator functions
\[ a(r,\rho, \eta) := \td a(r,[r]\rho, [r]\eta), \ 
   b(r, \rho, \eta) := \td b (r,[r]\rho, [r]\eta)   \]
will be interpreted as amplitude functions of a pseudo-differential calculus on $\R$ containing $\eta$ as a parameter (below we assume $\eta\not=0$). We intend to apply an analogue of Kumano-go's technique \cite{Kuma1} and form the oscillatory integral
\begin{equation} \label{cylnew9}
a \# b (r,\rho,\eta) =  \iint e^{-it \tau} a(r,
      \rho+\tau,\eta) b(r + t, \rho,\eta) dt \dbar \tau
\end{equation}
which has the meaning of a Leibniz product, associated with the composition of operators. The rule
\begin{equation} \label{cylnew29}
\op_r(a)(\eta)\op_r(b)(\eta)=\op_r(a\# b)(\eta)
\end{equation}
for $\eta\not=0$ will be justified afterwards.
Similarly as in \cite{Kuma1}, applying Taylor's formula, the function $a\# b$ can be decomposed in the form  
\begin{equation} 
\label{152.0809.eq}
 a \# b (r,\rho,\eta)  =  \ \sum_{k=0}^{N} \frac{1}{k!}
 \partial_\rho^k a(r,\rho,\eta) D_r^k b(r,\rho,\eta) +
 r_N(r,\rho,\eta) 
 \end{equation}
for
\begin{align}
\label{152.131neu2807.eq} 
 r_N(r,\rho,\eta)  = \ \frac{1}{N!} \iint e^{-it \tau} 
      \{\int_{0}^{1} (1-\theta)^N &(\partial_\rho^{N+1} a)
      (r,\rho+\theta \tau, \eta)d \theta \}  \\
    &  (D_r^{N+1} b)(r+t,\rho,\eta)dt \dbar \tau. \notag
\end{align}
By virtue of Corollary \ref{cylcormul} we have $\frac{1}{k!} 
\partial_\rho^k a(r, \rho,\eta) D_r^k
       b(r,\rho,\eta) =: c_k(r, \rho,\eta)$    	    
for $c_k (r, \rho,\eta) =$  $\td c_k(r,[r]\rho,[r]\eta)$,
$\td c_k(r, \td\rho,\td\eta) \in
S^{\nu+\td\nu} (\R,L_\clw^{\mu+\td\mu-k})$. Let us now
characterise the remainder. 

\begin{Lem}
\label{152.2407.le}
For every $ s', s''
\in \R$, $l,m,k \in \N$, there is an $N \in \N$ such 
that 
\begin{equation}
\label{152.(NE)1608.eq}
\| D_r^i D_\rho^j r_N(r,\rho,\eta) \|_{s',s''} \leq
    c \langle \rho \rangle^{-k} \langle r \rangle^{-l} 
    \langle \eta
    \rangle^{-m}
\end{equation}
for all $(r, \rho) \in \R^2$, $|\eta| \geq \varepsilon > 0$,
$i,j \in \N$, for some constant
 $c = c(s',s'',k,l,m,N, \varepsilon) > 0$, 
here $\| \cdot \|_{s',s''} = \| \cdot \|_{ {\Li}
(H^{s'}(X), H^{s''}(X))}$.    
\end{Lem}

\proof
Let us write $\boldsymbol{S}^{\mu,\nu} := \left\{ \td a 
(r,[r]\rho,[r]\eta) :
\td a(r,\td\rho,\td\eta) \in S^\nu(\R,
       L_\clw^\mu)\right\}$.    
By virtue of Proposition \ref{152.(cal)2209.pr} 
we have
\[ \partial_\rho^k \td a(r,[r]\rho,[r]\eta) 
    \in \boldsymbol{S}^{\mu-k, \nu+k},   
    \partial_r^k \td b(r,[r]\rho,[r]\eta)  \in
    \boldsymbol{S}^{\td\mu,\td\nu-k}    \]
for every $k$. Let us set
\begin{align*}
\td a_{N+1}(r,[r]\rho+[r]\theta \tau,[r]\eta) & :=    
 (\partial_\rho^{N+1} a) (r,\rho+\theta \tau,\eta),  \\
\td b_{N+1} (r+t,[r+t] \rho, [r+t]\eta) & := 
       (D_r^{N+1} b) (r+t,\rho,\eta).
\end{align*}        

By virtue of Theorem \ref{131.red2707.th} for every $s_0,s'' \in
\R$ and every $M$ there exists a $\mu(M)$ such that for every
$p(\td\rho,\td\eta) \in L_\clw^\mu(X;
\R_{\td\rho,\td\eta}^{1+q})$, $\mu \leq \mu(M)$, we
have
\begin{equation}
\label{152.133neu2209.eq}
  \| p(\td\rho, \td\eta) \|_{s_0,s''} \leq c
       \langle \td\rho, \td\eta \rangle^{- M}  
\end{equation}           
for all $(\td\rho,\td\eta) \in \R^{1+q}$, $c =
c(s_0,s'',\mu,M) > 0$. 
Moreover, for every $s',s_0 \in \R$ there exists a $B \in
\R$ such that $\|p(\td\rho,\td\eta)\|_{s',s_0}
\leq c \langle \td\rho, \td\eta \rangle^B$ for
all $(\td\rho, \td\eta) \in \R^{1+q}$, $c =
c(s',s_0,\mu) > 0$. 
We apply this for
$\td a_{N+1}(r,\td\rho, \td\eta)$ and
$\td b_{N+1}(r, \td\rho, \td\eta)$, combined with
the dependence on $r \in \R$ as a symbol in this variable. In other words,
we have the estimates

\begin{align}
\label{152.(E1)2807.eq}
\| \td a_{N+1}(r, \td\rho, \td\eta) \|_{s_0,s''}
     \leq c & \ \langle r \rangle^{\nu+(N+1)} \langle
     \td\rho,\td\eta \rangle^{- M},   \\
\label{152.(E2)2807.eq}
\| \td b_{N+1}(r,\td\rho,\td\eta) \|_{s',s_0}
     \leq & \ c \langle r \rangle^{\td\nu-(N+1)} \langle
     \td\rho, \td\eta \rangle^B;
\end{align}
here we applied the above-mentioned result to $\td a_{N+1}$
for the pair $(s_0,s'')$ for $N$ sufficiently large, and for
$\td b_{N+1}$ the second estimate for $(s',s_0)$ with some
exponent $B$. Let us take $s_0 := s'- \td\mu$; then we
can set $B = \max
\{\td\mu,0 \}$. The remainder \eqref{152.131neu2807.eq} is
regularised as an oscillatory integral in $(t,\tau)$, i.e., we may 
write 

\begin{align}
\label{152.rel2807.eq}
r_N& (r,\rho,\eta) =  \frac{1}{N!} \iint e^{-it \tau} \langle t
       \rangle^{-2L} (1- \partial_\tau^2)^L \langle \tau
       \rangle^{-2K} (1- \partial_t^2)^K   \\          
   & \ \Big\{ \int_{0}^{1} (1- \theta)^N \td a_{N+1}
        (r,[r]\rho +[r]\theta \tau,[r]\eta) d \theta\Big\}
	\td b_{N+1}(r+t,[r+t]\rho,[r+t]\eta) dt \dbar\tau
                            \notag
\end{align}
for sufficiently large $L,K$. For simplicity from now on we
assume $q = 1$; the considerations for the general case are
completely analogous. Then we have for every $l \leq L$
\[  \partial_\tau^{2l} \td a_{N+1} (r,[r]\rho+[r] 
      \theta \tau, [r]\eta) = \left(\partial_{\td\rho}^{2l}
      \td a_{N+1}\right) (r,[r]\rho+[r]\theta\tau,
      [r]\eta) ([r]\theta)^{2l},    \]
and for every $k \leq K$      
\begin{align*}
\partial_t^{2k} & \ \td b_{N+1} (r+t,[r+t]\rho, 
    [r+t]\eta) =
    \big(\partial_t^{2k} \td b_{N+1}\big) (r+t,[r+t]\rho,
    [r+t]\eta) \\
   & \ + \big(\partial_{\td\rho}^{2k} 
           \td b_{N+1}\big) 
       (r+t,[r+t]\rho,[r+t]\eta) (\rho 
       \partial_t[r+t])^{2k}
             \\
   & \ + \big(\partial_{\td\eta}^{2k} 
             \td b_{N+1}
	     \big) (r+t,[r+t]\rho, [r+t]\eta) 
	     (\eta \partial_t [r+t])^{2k} + R,
\end{align*}

where $R$ denotes several mixed derivatives.	     
From \eqref{152.(E1)2807.eq} we have
\begin{equation}
\label{152.(ES)2807.eq}
\| \partial_\tau^{2l} \td a_{N+1}
      (r,[r]\varrho+r[\theta]\tau,[r]\eta) \|_{s_0,s''}
      \leq c \langle r \rangle^{\nu+(N+1)} \langle [r]\varrho +
      [r]\theta \tau, [r]\eta \rangle^{-M-2l} ([r]\theta)^{2l},  
\end{equation}
see Corollary \ref{131.spec2707.co}, and \eqref{152.(E2)2807.eq}
gives us
\begin{equation}
\label{152.(ET)2807.eq}
\| (\partial_t^{2k} \td b_{N+1}) (r+t,[r+t] \varrho, [r+t]\eta) 
    \|_{s',s_0} \leq c \langle r+t \rangle^{\td\nu-(N+1)}
    \langle [r+t]\varrho, [r+t]\eta \rangle^B   
\end{equation}
(where we take $N$ so large that $\td\nu-(N+1) \leq 0$), and
\begin{align}
\label{152.(EU)2807.eq}
\| \bigl(\partial_{\td\varrho}^{2k} & \td b_{N+1} \bigr)
      (r+t,[r+t]\varrho, [r+t] \eta) (\varrho \partial_t [r+t])^{2k}
      \|_{s',s_0} \\
 & \leq c \langle r+t \rangle^{\td\nu-(N+1)} \langle [r+t]
      \varrho, [r+t]\eta \rangle^{B-2k} |\varrho 
      \partial_t[r+t]|^{2k},                 \notag
\end{align}
\begin{align}
\label{152.(EV)2807.eq}
\| \bigl( \partial_{\td\eta}^{2k} & \td b_{N+1}\bigr)
     (r+t,[r+t]\varrho,[r+t]\eta) (\eta \partial_t [r+t])^{2k}
     \|_{s',s_0}   \\
& \leq c \langle r+t \rangle^{\td\nu-(N+1)} \langle [r+t]\varrho,
[r+t]\eta \rangle^{B-2k} |\eta \partial_t [r+t]|^{2k}.   \notag
\end{align}
The above-mentioned mixed derivatives admit similar estimates (in fact,
better ones; so we concentrate on those contributed by
\eqref{152.(ES)2807.eq}, \eqref{152.(ET)2807.eq},
\eqref{152.(EU)2807.eq}, \eqref{152.(EV)2807.eq}).
    
We now derive an estimate for $\|r_N (r,\varrho,\eta) \|_{s',s''}$.
Using the relation \eqref{152.rel2807.eq} we have 
\begin{align*}
\| r_N& (r,\varrho,\eta) \|_{s',s''} \leq  \ \iint \int_{0}^{1} \|
    \langle t \rangle^{-2L} (1-\partial_\tau^2)^L \langle \tau
    \rangle^{-2K} (1- \partial_t^2)^K    \\
& \ (1-\theta)^N \td a_{N+1} (r,[r]\varrho+[r]\theta \tau, 
    [r]\eta)
       \td b_{N+1} (r+t,[r+t]\varrho, [r+t]\eta) \|_{s',s''} d
       \theta dt \dbar \tau.
\end{align*}
The operator norm under the integral can be estimated by expressions of
the kind  
\begin{align*}
I := & \ c \langle r \rangle^{\nu+(N+1)} \langle r+t
      \rangle^{\td\nu-(N+1)} \langle t \rangle^{-2L}
      \langle \tau \rangle^{-2K} \langle
      [r]\rho+[r]\theta \tau, [r]\eta 
      \rangle^{-M-2l} ([r]\theta)^{2l}   \\
  & \langle [r+t] \rho, [r+t] \eta \rangle^B \big\{ 1+
      \langle [r+t] \rho, [r+t] \eta \rangle^{-2k} (|
      \rho|^{2k}+ |\eta |^{2k})
          |(\partial_t [r+t])^{2k}| \big\}
\end{align*}
$l \leq L$, $k \leq K$, plus terms from $R$ of a similar
character.	      
We have, using Peetre's inequality,
\[ \langle r \rangle^{\nu+(N+1)} \langle r+t
     \rangle^{\td\nu-(N+1)} \leq \langle r 
     \rangle^{\nu+\td\nu} \langle t
     \rangle^{|\td\nu-(N+1)|}.    \]
Moreover, we have 
$\langle [r] \rho+[r]\theta \tau, [r] \eta \rangle^{-2l} 
      ([r]\theta)^{2l} 
      \leq c \langle [r] \eta \rangle^{-2l} [r]^{2l} \leq c$  
for $|\eta| \geq \varepsilon > 0$ (as always, $c$ denotes different
constants), and
\begin{multline*}
\langle [r+t]\rho, [r+t]\eta \rangle^{-2k}
    (|\rho|^{2k}+|\eta|^{2k}) |(\partial_t[r+t])^{2k}|  \\ 
 \leq  \ c \big\{ \langle [r+t]\rho \rangle^{-2k} 
         ([r+t]|\rho|)^{2k}
         + \langle[r+t]\eta \rangle^{-2k} ([r+t]|\eta|)^{2k} \big\}
	   [r+t]^{-2k} \leq c,
\end{multline*}
using $|(\partial_t[r+t])^{2k}| \leq c, [r+t]^{-2k} \leq c$ for all
$r,t \in \R$ and $|\zeta| \leq c \langle \zeta \rangle$ for every
$\zeta$ in $\R^d$.
This yields
\[ I \leq  c \langle r \rangle^{\nu+\td\nu} \langle t
         \rangle^{|\td\nu-(N+1)|} \langle t \rangle^{-2L}
	 \langle \tau \rangle^{-2K} 
    \langle [r]\rho + [r] \theta \tau, [r] \eta \rangle^{- M}
       \langle[r+t]\rho, [r+t]\eta \rangle^B.\]	    	   
Writing $M = M' + M''$ for suitable $M',M'' \geq 0$ to be fixed later on, we have 
\begin{align*}
\langle  [r] \rho & + [r] \theta \tau, [r] \eta \rangle^{-M} = 
      \langle [r] \rho + [r]\theta \tau, [r]\eta \rangle^{-M'}
         \langle [r] \rho + [r] \theta \tau, [r] \eta
	 \rangle^{-M''}                \\
 \leq & c \langle  [r] \eta \rangle^{-M'} \langle [r] \rho, [r] 
 \eta
          \rangle^{-M''} \langle [r] \theta \tau \rangle^{M''}     
 \leq c \langle [r] \eta \rangle^{- M'} \langle [r] \rho
          \rangle^{-M''} \langle [r] \theta \tau 
	     \rangle^{M''}.
\end{align*}	  	  
We applied once again Peetre's inequality which gives us also
\[ \langle [r+t]\rho, [r+t] \eta \rangle^B \leq c \langle [r+t]
       \rho \rangle^B \langle [r+t] \eta \rangle^B     \]	  
since $B \geq 0$. Thus
\begin{multline*}
I \leq \ c \langle r \rangle^{\nu+\td\nu} \langle t
     \rangle^{|\td\nu-(N+1)|-2L}\langle \tau \rangle^{-2 K}\langle [r] \theta \tau \rangle^{M''} \\
   \ \langle [r+t] \rho \rangle^B \langle [r] \rho 
       \rangle^{- M''} \langle [r+t] \eta \rangle^B \langle [r]\eta\rangle^{-M'}.
\end{multline*}

Let us show 
\[  \langle t \rangle^{- B} \langle [r+t]   \rho \rangle^B
        \langle [r] \rho \rangle^{- B} \leq c.     \]
In fact, this is evident in the regions $|r| \leq C$, $|t| \leq C$
or $|r| \geq C$, $|t| \leq C$ for some $C > 0$. For $|r| \leq C$,
$|t| \geq C$ the estimate essentially follows from $1+t^2 \rho^2
\leq (1+t^2)(1+ \rho^2)$. For $|r| \geq C$, $|t| \geq C$, $[r+t]
\leq C$ the estimate is evident as well. It remains the case $|r|
\geq C$, $|t| \geq C$, $[r+t] \geq C$, where the estimate follows (for $C \geq 1$
so large that $[r+t] = |r+t|$, $[r] = |r|$) from 
\begin{align*}
\langle t \rangle^{-2} & \langle [r+t]\rho \rangle^2 \langle [r]
    \rho \rangle^{-2}   
 =  \frac{1+|r+t|^2 | \rho|^2}{(1+ |t|^2)(1+|r \rho|^2)}
     \leq  \frac{1+|r \rho|^2 + 2|rt \rho|^2 |t \rho|^2} {1+
     |t|^2 + |r \rho|^2 + |rt \rho|^2}   \\
 \leq & \ c \frac{1+ |r \rho |^2 + |t \rho |^2 + 2 | r  t
          \rho|^2}  {1+ |r \rho |^2 + | t \rho |^2}
	\leq c \Bigl( 1+ \frac{2 | rt \rho|^2}  {1+ |r \rho|^2
	   + |t \rho|^2} \Bigr)  \leq \textup{const}.
\end{align*}
Here we employed $|r t \rho |^2 \geq |t \rho |^2$ for $|r|
\geq C \geq 1$  and
\[ \frac{|rt\rho|^2} {1+|r\rho|^2 + |t \rho|^2}   \leq 
       \frac{r^2 t^2} { r^2 + t^2}  =  \frac{r^2} {r^2 + t^2}
         \frac{t^2}{r^2 + t^2}  \leq \textup{const}.     \]
Analogously we have $\langle t \rangle^{- B} \langle [r+t] \eta 
\rangle^B \langle [r]
     \eta \rangle^{- B} \leq c$.   
This gives us the estimate
\[   I \leq c \langle r \rangle^{\nu+\td\nu} \langle t
\rangle^{|\td\nu-(N+1)| - 2L + 2B} \langle\tau \rangle^{-2K}\langle [r]\theta \tau \rangle^{M''} \langle [r]\rho\rangle^{B-M''} \langle [r]\eta \rangle^{B-M'}  .     \]
Finally, using $\langle \tau \rangle^{-M''} \langle r 
\rangle^{-M''} \langle[r]\theta \tau \rangle^{M''} \leq c$    	     	 
for all $0 \leq \theta \leq 1$ and all $r, \tau$, we obtain 
\[  I \leq c \langle r \rangle^{\nu+\td\nu+M''} \langle t
        \rangle^{|\td\nu-(N+1)|-2L + 2B} \langle \tau
	\rangle^{-2K+M''} \langle [r] \rho \rangle^{B-M''} \langle
	[r] \eta \rangle^{B-M'}     \]	 	         
for all $r,t \in \R$, $\rho, \tau \in \R$, $0 \leq \theta \leq
1$. Choosing $K$ and $L$ so large that
\[  -2K+M'' < - 1, |\td\nu-(N+1)| + 2B - 2 L < - 1,    \]
it follows that 
$\| r_N(r, \rho,\eta) \|_{s',s''} \leq c \langle r
       \rangle^{\nu+\td\nu+M''} \langle [r]\eta\rangle^{B-M'}
       \langle \rho\rangle^{B-M''}$ for $\eta \not= 0$ using that $\langle 
[r] \rho \rangle^{B-M''} \leq c \langle
\rho \rangle^{B-M''}$ for $B-M'' \leq 0$.
Let us now show that for $B-M'\leq 0$
\begin{equation}
\label{152.(est)1308.eq}
\langle [r]\eta \rangle^{B-M'} \leq c [r]^{B-M'} \langle \eta
       \rangle^{B-M'}
\end{equation}
for all $|\eta| \geq \varepsilon > 0$ and some $c = c (\varepsilon)
> 0$. In fact, we have
\[
\frac{[r]^2\langle \eta \rangle^2}{1+|[r]\eta|^2} =
     \frac{[r]^2}{1+|[r]\eta|^2}
       \frac{\langle \eta \rangle^2}{ 1+|[r]\eta|^2}  \leq c
          \frac{1}{[r]^{-2}+|\eta|^2}   
	    \frac{1}{|\eta|^{-2}+[r]^{-2}} \leq c,  \]
i.e., $(1+|[r]\eta|^2)^{-1} \leq c [r]^{-2} \langle \eta
\rangle^{-2}$ which entails the estimate
\eqref{152.(est)1308.eq}.
It follows
\[ \|r_N(r,\rho,\eta)\|_{s',s'}\leq c\la r\ra^{\nu+\td\nu+M''+B-M'}\la\rho\ra^{B-M''} \la\eta\ra^{B-M'} \] 
Now $B$ is fixed, and $M,M''$ can be chosen independently so large that 
\[ B-M''\leq -k, \ B-M'\leq -m, \ \nu+\td\nu+M''+B-M'\leq -l. \] 
Therefore, we proved that for every $s',s'' \in \R$ and $k,l,m \in \N$ there is an $N \in \N$ such that
\begin{equation}
\label{152.135neu11308.eq}  
\| r_N(r, \rho,\eta) \|_{s',s''} \leq c \langle \rho
       \rangle^{-k} \langle r \rangle^{-l} \langle \eta \rangle^{-m}     
\end{equation}
for all $(r, \rho) \in \R^2$, $|\eta| \geq \varepsilon > 0$. 
In an analogous manner we can show the estimates 
\eqref{152.(NE)1608.eq} for all $i,j$. 
\qed            
\begin{Rem} \label{cylremmdeg}
For future references we call an operator $\op_r(r_N)(\eta)$ for an operator function $r_N(r,\rho,\eta)$ satisfying the estimate \eqref{152.135neu11308.eq} smoothing of degree $(k,l)$ \lr with respect to given fixed $s',s''$\rr. A similar notation makes sense and will be used below when we ignore $\eta$ and replace $\la\eta\ra^{-m}$ by some fixed constant $>0$.
\end{Rem}
\begin{Prop}
\label{152.2209.co}
For every $a(r,\rho,\eta) \in S^{\mu,\nu}$ and $\varphi(r) =
[r]^{\td\nu}$ \lr which belongs to 
$\boldsymbol{S}^{0,{\td\nu}}$\rr\ for every $\eta \not= 0$ we have 
\lr as operators $\op_r\big(\td a(r,[r]\rho,[r]\eta)\big):C_0^\infty\big(\R,C^\infty(X)\big)\to C^\infty\big(\R,C^\infty(X)\big)$\rr
\begin{equation} \label{cylnew13}  
\op_r(a)(\eta) \circ \varphi = \varphi \circ \op_r(b)(\eta) + d(\eta)
\end{equation}
for some $b(r, \rho, \eta) \in \boldsymbol{S}^{\mu, \nu}$ and a
remainder $d(\eta) = \op_r(r_N)(\eta)$ which is an operator 
function $r_N(r, \rho,
\eta) \in C^\infty\big(\R \times \R \times \R_\eta^q$, $\Li\big(H^{s'}(X), H^{s''}(X)\big)\big)$ for every given $s', s''$ and
sufficiently large $N = N(s',s'') \in \N$,  satisfying the
estimates \eqref{152.(NE)1608.eq} for all $(r,\rho) \in
\R^2$ and all $|\eta|\geq \varepsilon > 0$.     
\end{Prop}
\proof
We apply the relation \eqref{152.0809.eq} to the case $b(r,
\rho,\eta) = \varphi(r)$ and obtain
\[  \op(a) \circ \varphi = \op(a \# \varphi) = \sum_{k=0}^{N}
      \op \big( \frac{1}{k!} \partial_\rho^k a(r,\rho,\eta)
      D_r^k \varphi(r)\big) + \op(r_N).    \]
According to Corollary \ref{cylcormul} we can form $\big(
\partial_\rho^k a D_r^k \varphi\big)
(r,[r]\rho,[r]\eta)$ with
$\big( \partial_\rho^k a D_r \varphi\big)(r,
\td\rho,\td\eta) \in S^{\nu+\td\nu}(\R,
L_\clw^{\mu-k}(X; \R^{1+q}))$. 
There is then a $\td c_N(r,\td\rho,\td\eta)
\in S^{\nu+\td\nu} (\R, L_\clw^{\mu-(N+1)}(X; \R^{1+q}))$
which is the asymptotic sum of the symbols
$\frac{1}{k!}\big(\partial_\rho^k a D_r^k \varphi\big) (r,
\td\rho,\td\eta)$ over $k\geq N+1$. Writing
$c_N(r,\rho,\eta) = \td c_N(r,[r]\rho,[r]\eta)$ we
obtain
\[  a \# b (r,\rho,\eta) = p_N(r,\rho,\eta) +
       d_N(r,\rho,\eta)    \]
for $p_N(r,\rho,\eta) = \td p_N(r,[r]\rho,[r]\eta)$, 
\[  p_N(r, \rho,\eta) = \sum_{k=0}^{N} \frac{1}{k!}
      \partial_\rho^k a(r,\rho, \eta) D_r^k \varphi(r) + 
      c_N (r,\rho,\eta) \in \boldsymbol{S}^{\mu,\nu+\td\nu}  \] 
and $d_N(r, \rho,\eta) = r_N(r, \rho,\eta) -
c_N(r,\rho,\eta)$. Now $r_N(r,\rho,\eta)$ satisfies the 
desired estimates. 
Similarly as in connection with \eqref{152.133neu2209.eq} for
every $s',s'' \in \R$ and $M \in \N$ we find an $N\in\N$ sufficiently
large such that
\[  \| c_N(r,\td\rho,\td\eta) \|_{s',s''} \leq c
       \la r \ra^{\nu+\td\nu} \la 
       \td\rho, \td\eta \ra^{- 4 M}.      \]
This entails.
\[  \| c_N(r,[r]\rho, [r]\eta) \|_{s',s''} \leq c \la r
    \ra^{\nu+\td\nu} \la [r]\rho,
    r[\eta]\ra^{-4 M}     \]
for all $r,\rho,\eta$. Now
\begin{align*}  \la [r]\rho, r [\eta] \ra^{-4} = & \ 
               [r]^{-4}
     \bigg( \frac{1}{\frac{1+[r]^2 \rho^2 + [r]^2\eta^2}
         {[r]^2}}\bigg)^2    \\
     = & \ [r]^{-4} \frac{1}{[r]^{-2} + \rho^2 + \eta^2}
          \frac{1}{[r]^{-2} + \rho^2 + \eta^2}
     \leq c [r]^{-4} \la \rho \ra^{-2} \la \eta
     \ra^{-2}    
\end{align*}	           
for $|\eta| \geq \varepsilon > 0$, for a constant $c =
c(\varepsilon ) > 0$. We thus obtain
\[  \|c_N(r,[r]\rho, [r]\eta) \|_{s',s''} \leq c \la r
      \ra^{\nu+\td\nu-4M} \la \rho \ra^{-2M}
      \la \eta \ra^{-2 M}.    \]
This completes the proof since $M$ is arbitrary.                       
\qed \\
\indent Let us now return to the interpretation of \eqref{cylnew9} as the left symbol of a composition of operators. From Theorem \ref{152.2309.th} we know that
\[\op_r(a)(\eta),\op_r(b)(\eta):\Sch\big(\R,C^\infty(X)\big)\to \Sch\big(\R,C^\infty(X)\big) \]
are continuous operators. Thus also $\op_r(a)(\eta)\op_r(b)(\eta)$ is continuous between the Schwartz spaces. This shows, in particular, that the oscillatory integral techniques of \cite{Kuma1} also apply for our (here operator-valued) amplitude functions, and we obtain the relation \eqref{cylnew29}.\\
\indent Let $A(\eta)=\op_r(a)(\eta)$ for 
\[ a(r,\rho,\eta):=\td a(r,[r]\rho,[r]\eta),\ \td a(r,\td\rho,\td\eta)\in \boldsymbol{S}^{\mu,\nu}. \]
Then we form the formal adjoint $A^*(\eta)$ with respect to the $L^2(\R\times X)$-scalar product, according to
\[ \big(A(\eta)u,v\big)_{L^2(\R\times X)}=\big(u,A^*(\eta)v\big)_{L^2(\R\times X)} \]
for all $u,v\in\Sch\big(\R,C^\infty(X)\big)$. As usual we obtain
\[ A^*(\eta)v(r')=\op_{r'}(a^*)(\eta)v(r') \]
for the right symbol $a^*(r',\rho,\eta)=\bar a(r',\rho,\eta)=\td{\bar a}(r',[r']\rho,[r']\eta)$. Similarly as before we can prove that 
\[ \op_{r'}(a^*)(\eta):\Sch\big(\R,C^\infty(X)\big)\to\Sch\big(\R,C^\infty(X)\big) \]
is continuous for every $\eta\not=0$. Thus by duality it follows that
\begin{equation} \label{cyluen15} 
\op_r(a)(\eta):\Sch'\big(\R,\E'(X)\big)\to\Sch'\big(\R,\E'(X)\big)
\end{equation}
is continuous for every $\eta\not=0$. Note here that $f\in\E'(X)\Leftrightarrow f\in H^s(X)$ for some real $s\in\R$; then $\Sch'\big(\R,\E'(X)\big)$ means the inductive limit of the spaces $\Li\big(\Sch(\R),H^s(X)\big)$ over $s\in\R$.
\begin{Lem} \label{cyllemrest}
For every $s',s''\in\R$ and $l,m,k\in\N$ there exists a real $\mu(s',s'',k,l,m)$ such that for every $a(r,\rho,\eta)\in\boldsymbol{S}^{\mu,0}$ we have
\[ \|a(r,\rho,\eta)\|_{s',s''}\leq c \la\rho\ra^{-k}\la r\ra^{-l}\la\eta\ra^{-m} \]
whenever $\mu\leq\mu(s',s'',k,l,m)$, $|\eta|\geq\varepsilon>0$.
\end{Lem}
\proof
The proof is straightforward, using Theorem \ref{131.red2707.th}, more precisely, writing $a(r,\rho,\eta)=\td a(r,[r]\rho,[r]\eta)$, we have the estimate
\[\|\td a(r,\td\rho,\td\eta)\|_{s',s''}\leq c\la\td\rho,\td\eta\ra^{-N} \]
for every fixed $N\in\N$ when $\mu$ is chosen sufficiently negative (depending on $N$), uniformly in $r\in\R$. Then, similarly as in the proof of Lemma \ref{152.2407.le}, we obtain for suitable $N$ and given $k,l,m$ that $\la[r]\rho,[r]\eta\ra^{-N}\leq c\la\rho\ra^{-k}\la r\ra^{-l}\la\eta\ra^{-m}$ for $|\eta|\geq\varepsilon>0$.
\qed
\begin{Thm} \label{cylCalVai}
For every $\td p(r,\td\rho,\td\eta)\in S^0\big(\R,L_\clw^s(X;\R^{1+q}_{\td\rho,\td\eta})\big)$, $s\leq 0$, and $p(r,\rho,\eta) = \td p(r,[r]\rho,[r]\eta)$, the operator
\begin{equation} \label{cylcont1}
\op_r(p)(\eta):L^2(\R\times X)\to L^2(\R\times X)
\end{equation}
is continuous for every $\eta\in\R^q\setminus\{0\}$, and we have
\begin{equation} \label{cylop1}
\|\op_r(p)(\eta)\|_{\Li(L^2(\R\times X))}\leq c\la\eta\ra^s
\end{equation}
for all $|\eta|\geq\varepsilon, \varepsilon >0$ and a constant $c=c(\varepsilon)>0$.
\end{Thm}
\proof
For the continuity \eqref{cylcont1} and the estimate \eqref{cylop1} we apply a version of the Calder\'on-Vaillancourt theorem which states that if $H$ is a Hilbert space and $a(r,\rho)\in C^\infty\big(\R\times\R,\Li(H)\big)$ a symbol satisfying the estimate
\begin{equation} \label{cylnorm1}
\pi(a):=\sup_{\substack{k,l=0,1 \\ (r,\rho)\in\R^2}}\|D_r^lD_\rho^k a(r,\rho)\|_{\Li(H)}<\infty
\end{equation}
the operator
\[ \op_r(a):L^2(\R,H)\to L^2(\R,H) \]
is continuous, where
\[ \|\op_r(a)\|_{\Li(L^2(\R,H))}\leq c\pi(a) \]
for a constant $c>0$. In the present case we have
\begin{equation} \label{cyleqp1}
a(r,\rho)=p(r,[r]\rho,[r]\eta)
\end{equation}
where $\eta\not=0$ appears as an extra parameter. It is evident that the right hand side of \eqref{cyleqp1} belongs to $C^\infty\big(\R\times\R\times\R^q,\Li(L^2(X))\big)$. From the assumption on $\td p(r,\td\rho,\td\eta)$ we have 
\begin{equation} \label{cylep1}
\sup_{r\in\R}\|\td p(r,\td\rho,\td\eta)\|_{\Li(L^2(X))}\leq c\la\td\rho,\td\eta\ra^s
\end{equation}
for all $(\td\rho,\td\eta)\in\R^{1+q}$ and some $c>0$. In fact, when $\td p$ is independent of $r$ the latter estimate corresponds to \eqref{cyleqlam} for $s=\nu=0$ and $\mu=s\leq 0$. In the $r$-dependent case the operator norms that play a role in Theorem \ref{131.red2707.th} are uniformly bounded in $r\in\R$, since $\td p(r,\td\rho,\td\eta)$ is a symbol of order 0 in $r$ with values in $L_\clw^s(X;\R^{1+q}_{\td\rho,\td\eta})$. 
For \eqref{cylnorm1} we first check the case $l=k=0$.
We have
\begin{equation} \label{cylesti1}
\sup_{(r,\rho)\in\R^2}\la[r]\rho,[r]\eta\ra^s\leq c\la\eta\ra^s
\end{equation}
for all $|\eta|\geq\varepsilon>0$ and some $c=c(\varepsilon)>0$. Thus \eqref{cylep1} gives us
\[ \sup_{(r,\rho)\in\R^2}\|\td p(r,[r]\rho,[r]\eta)\|_{\Li(L^2(X))} \leq c\la\eta\ra^s \]
for such a $c(\varepsilon)>0$. Assume now for simplicity $q=1$ (The general case is analogue). For the first order derivatives of $\td p(r,[r]\rho,[r]\eta)$ in $r$ we have 
\begin{equation} \label{cyleqsu1}
\partial_r \td p(r,[r]\rho,[r]\eta)=(\partial_r\td p)(r,[r]\rho,[r]\eta)+[r]'(\rho\partial_{\td\rho}+\eta\partial_{\td\eta})\td p(r,[r]\rho,[r]\eta) 
\end{equation}
for $[r]'=\frac{d}{dr}[r]$. For the derivatives of $\td p$ with respect to $\td\rho,\td\eta$ we employ that $\partial_{\td\rho}\td p(r,\td\rho,\td\eta),\partial_{\td\eta}\td p(r,\td\rho,\td\eta)\in S^0\big(\R,L^{s-1}(X;\R^{1+q}_{\td\rho,\td\eta})\big)$. Thus, similarly as before we obtain 
\[ \|\partial_{\td\rho,\td\eta}^\alpha \td p(r,\td\rho,\td\eta)\|_{\Li(L^2(X))}\leq c\la\td\rho,\td\eta\ra^{s-1} \]
for any $\alpha\in\N^2$, $|\alpha|=1$. This gives us for the summand on the right of \eqref{cyleqsu1}
\begin{multline*}
\sup_{(r,\rho)\in\R^2}\|[r]^{-1}[r]'([r]\rho\partial_{\td\rho}+[r]\eta\partial_{\td\eta})p(r,[r]\rho,[r]\eta)\|_{\Li(L^2(X))} \\ \leq \sup [r]^{-1}|[r]\rho+[r]\eta|\la r\rho,r\eta\ra^{s-1} \\ \leq c\la\eta\ra^s\sup [r]^{-1}|[r]\rho,[r]\eta|\la r\rho,r\eta\ra^{-1} \leq c\la\eta\ra^s.
\end{multline*}
Here we employed \eqref{cylesti1}. For the derivative of $p(r,[r]\rho,[r]\eta)$ in $\rho$ we have
\begin{multline*} 
\sup\|\partial_\rho \td p(r,[r]\rho,[r]\eta)\|_{\Li(L^2(X))}=\sup \|[r](\partial_{\td\rho}\td p)(r,[r]\rho,[r]\eta)\|_{\Li(L^2(X))} \\ \leq c\sup[r]\la[r]\rho,[r]\eta\ra^{s-1} \leq c\la\eta\ra^s
\end{multline*}
for all $|\eta|\geq\varepsilon > 0$. This gives altogether the estimate \eqref{cylop1}.
\qed
\begin{Thm} \label{cylthminj}
Let $\td p(\td\rho,\td\eta)\in L_\clw^s(X;\R^{1+q})$ be parameter-dependent elliptic of order $s\in\R$, and set $p(r,\rho,\eta)=\td p([r]\rho,[r]\eta)$. Then there exists a $C>0$ such that for every $|\eta|\geq C$ the operator
\begin{equation} \label{cyleqS}
[r]^{-s}\op_r(p)(\eta):\Sch\big(\R,C^\infty(X)\big)\to\Sch\big(\R,C^\infty(X)\big)
\end{equation}
extends to an injective operator
\begin{equation} \label{cyleqT}
[r]^{-s}\op_r(p)(\eta):L^2(\R\times X)\to\Sch'\big(\R,\E'(X)\big).
\end{equation}
More precisely, considering $[r]^{-s}\op_r(p)(\eta)$ as an operator
\begin{equation} \label{cyleqCON}
[r]^{-s}\op_r(p)(\eta):L^2(\R\times X)\to\Li\big(\la r\ra^{-g}H^l(\R),H^t(X)\big)
\end{equation}
which is continuous for some $t\in\R$ and all $g,l\in\R$, then it is injective.
\end{Thm}
\proof
First, according to \eqref{cyluen15} there is a $t$ such that \eqref{cyleqCON} is continuous for all $g,l\in\R$. For the injectivity we show that the operator hast a left inverse. This will be approximated by $\op_r(a)$ for 
\begin{equation} \label{cylqea1}
a(r,\rho,\eta):=[r]^s\td p^{(-1)}([r]\rho,[r]\eta)
\end{equation}
where $\td p^{(-1)}(\td\rho,\td\eta)\in L_\clw^{-s}(X;\R^{1+q})$ is a parameter-dependent parametrix of $\td p(\td\rho,\td\eta)$. Setting
\begin{equation} \label{cylqeb1}
b(r,\rho,\eta):=[r]^{-s}\td p([r]\rho,[r]\eta)
\end{equation}
we can write the composition of the associated pseudo-differential operators in $r$ for every $N\in\N$ we have in the form 
\begin{equation} \label{cylqeLei}
\op_r(a)(\eta)\op_r(b)(\eta)=\op_r(a\# b)(\eta)=\op_r(1+c_N(r,\rho,\eta)+r_N(r,\rho,\eta))
\end{equation}
for $c_N(r,\rho,\eta)=\sum_{k=1}^N\frac{1}{k!}\partial_\rho^k a(r,\rho,\eta)D_r^k b(r,\rho,\eta)$ which has the form $c_N(r,\rho,\eta)=\td c_N(r,[r]\rho,[r]\eta)$ for some $\td c_N(r,\td\rho,\td\eta) \in S^0\big(\R,L_\clw^{-1}(X;\R^{1+q})\big)$. Moreover, the remainder $r_N$ is as in \eqref{152.131neu2807.eq}. From Theorem \ref{cylCalVai} for $s=-1$ we know that
\[ \|\op_r(c_N)(\eta)\|_{\Li(L^2(\R\times X))}\leq c\la\eta\ra^{-1} \]
for $|\eta|>\varepsilon$. Moreover, Lemma \ref{152.2407.le}, applied to $s'=s''=0$ together with an operator-valued version of the Calder\'on-Vaillancourt theorem, gives us 
\[ \|\op_r(r_N)(\eta)\|_{\Li(L^2(\R\times X))}\leq c\la\eta\ra^{-1} \]
for sufficiently large $N$. Thus for every $|\eta|$ sufficiently large the operator on the right of \eqref{cylqeLei} is invertible in $L^2(\R\times X)$, i.e., $\op_r(b)(\eta)$ has a left inverse which implies the injectivity. \qed
\begin{Rem}
Theorem \textup{\ref{cylthminj}} can easily be generalised to $\td p(r,\td\rho,\td\eta)\in S^0\big(\R,$ $L_\clw^s (X;\R^{1+q})\big)$, parameter-dependent elliptic for every fixed $r\in\R$, such that there exists a $\td p^{(-1)}(r,\td\rho,\td\eta)\in S^0\big(\R,L_\clw^{-s}(X;\R^{1+q})\big)$ which is a parameter-dependent parametrix for every fixed $r\in\R$. In other words, for $p(r,\rho,\eta)=\td p(r,[r]\rho,[r]\eta)$ there is a $C>0$ such that \eqref{cyleqS} extends to an injective operator \eqref{cyleqT} for every $|\eta|\geq C$.
\end{Rem}
In fact, if we form \eqref{cylqea1} and \eqref{cylqeb1} in an analogous manner including the extra $r$-dependence, we may apply the relation \eqref{cylqeLei}. The remaining conclusions in the proof of Theorem \eqref{cylthminj} do not depend on the assumption that $\td p$ is independent of $r$.

\section[Operators on an infinite cylinder]{Parameter-dependent operators on an infinite cylinder}

\subsection{Weighted cylindrical spaces}
\begin{Def} \label{cylDEFH}
Let $\td p(\td\rho,\td\eta)\in L_\clw^s(X;\R_{\td\rho,\td\eta}^{1+q})$ be as in Theorem \textup{\ref{cylthminj}}. Then $H_\textup{cone}^{s,g}(X^\asymp)$ for $s,g\in\R$ is defined to be the completion of $\Sch(\R\times X)$ with respect to the norm
\[ \|[r]^{-s+g}\op_r(p)(\eta^1)u\|_{\Li(L^2(\R\times X))} \]
for any fixed $\eta^1\in\R^q$, $|\eta^1|\geq C$ for some $C>0$ sufficiently large.
\end{Def}
\indent Setting $p^{s,g}(r,\rho,\eta):=[r]^{-s+g}\td p([r]\rho,[r]\eta)$, from Definition \ref{cylDEFH} it follows that 
\[ \op(p^{s,g})(\eta^1):\Sch(\R\times X)\to\Sch(\R\times X) \]
extends to a continuous operator
\begin{equation} \label{cylqeeord} 
\op(p^{s,g})(\eta^1):H^{s,g}_\textup{cone}(X^\asymp)\to L^2(\R\times X).
\end{equation}
\begin{Thm} \label{cylthmiso}
The operator \eqref{cylqeeord} is an isomorphism for every fixed $s,g\in\R$ and $|\eta^1|$ sufficiently large.
\end{Thm}
\proof
We show the invertibility by verifying that there is a right and a left inverse. By notation we have $p^{s,g}(r,\rho,\eta)=[r]^{-s+g}\td p([r]\rho,[r]\eta)\in\boldsymbol{S}^{s,-s+g}$. The operator family $\td p(\td\rho,\td\eta)\in L_\clw^s(X;\R^{1+q}_{\td\rho,\td\eta})$ is invertible for large $|\td\rho,\td\eta|\geq C$ for some $C>0$. There exists a parameter-dependent parametrix $\td p^{(-1)}(\td\rho,\td\eta)\in L_\clw^{-s}(X;\R^{1+q}_{\td\rho,\td\eta})$ such that $\td p^{(-1)}(\td\rho,\td\eta)=\td p^{-1} (\td\rho,\td\eta)$ for $|\td\rho,\td\eta|\geq C$. Let us set
\[ p^{-s,-g}(r,\rho,\eta):=[r]^{s-g}\td p^{(-1)}([r]\rho,[r]\eta)\in \boldsymbol{S}^{-s,s-g}, \]
and $P^{s,g}(\eta):=\op(p^{s,g})(\eta)$, $P^{-s,-g}(\eta):=\op(p^{-s,-g})(\eta)$. Then we have 
\begin{equation} \label{cyleeqinv}
P^{s,g}(\eta)P^{-s,-g}(\eta)=1+\op(c_N)(\eta)+R_N(\eta)
\end{equation}
for some $c_N(r,\rho,\eta)\in\boldsymbol{S}^{-1,0}$ and a remainder $R_N(\eta)=\op(r_N)(\eta)$ where $r_N$ is as in Lemma \ref{152.2407.le}. We have $\op(c_N)(\eta)\to 0$ and $R_N(\eta)\to 0$ in $\Li(L^2(\R\times X))$ as $|\eta|\to\infty$; the first property is a consequence of Theorem \ref{cylCalVai}, the second one of the estimate \eqref{152.(NE)1608.eq}. Thus \eqref{cyleeqinv} shows that $P^{s,g}(\eta)$ has a right inverse for $|\eta|$ sufficiently large. Such considerations remain true when we interchange the role of $s,g$ and $-s,-g$. In other words, we also have
\[ P^{-s,-g}(\eta)P^{s,g}(\eta)=1+\op(\td c_N)(\eta)+\widetilde R_N(\eta) \]
where $\op(\td c_N)(\eta)$ and $\widetilde R_N(\eta)$ are of analogous behaviour as before. This shows that $P^{s,g}(\eta)$ has a left inverse for large $|\eta|$, and we obtain altogether that \eqref{cylqeeord} is an isomorphism for $\eta=\eta^1$, $|\eta^1|$ sufficiently large.
\qed

\subsection{Elements of the calculus}
The results of Section \ref{cylsecleib} show the behaviour of compositions of parameter-dependent families $\op(a)(\eta)$ for $a(r,\rho,\eta)\in\boldsymbol{S}^{\mu,\nu}$ and $\eta\not=0$, first on $\Sch(\R\times X)$. In particular, Lemma \ref{152.2407.le} suggests the nature of smoothing symbols in such a calculus, namely, to be functions in $\eta\not=0$ with values in $\Sch\big(\R_r \times \R_\rho,\Li\big(H^{s'}(X),H^{s''}(X)\big)\big)$. Another interpretation of the results is that we can open an operator algebra on $X^\asymp$ either as an algebra of $\eta$-dependent families or where the value of $\eta\not=0$ is fixed in different ways depending on the operator. In the latter case Lemma \ref{152.2407.le} gives us remainders $K$ where $\|\cdot\|_{s',s''}$ is only estimated by $c\la\rho\ra^{-k}\la r\ra^{-l}$ without an explicit presence of $\eta$. It can be proved that, when we concentrate, for instance, on the case $s'=s''=0$, invertible operators of the form $1+K:L^2(\R\times X)\to L^2(\R\times X)$ can be written in the form $1+L$ where $L$ is again an operator of such a behaviour. Moreover, the composition of such a (smoothing) operator with an operator $\op(a)(\eta)$, $a\in\boldsymbol{S}^{\mu,\nu}$, $\eta\not=0$ fixed, gives us again an operator, smoothing in that sense. This can easily be deduced from the estimate \eqref{cyleqlam} of Theorem \ref{131.red2707.th}. Moreover, there are other (more or less standard) constructions that are immediate by the results of Section \ref{cylsec1}. For instance, if we look at $c(r,\rho,\eta)\in\boldsymbol{S}^{-1,0}$ in the relation \eqref{cyleeqinv}, by a formal Neumann series argument we find a $d(r,\rho,\eta)\in\boldsymbol{S}^{-1,0}$ such that 
\[ \big(1+\op(c)\big)\big(1+\op(d)\big)=1+\op(r_N) \]
for every $N\in\N$ with a remainder $r_N$ which is again as in Lemma \ref{152.2407.le}.
\begin{Thm}
Let $a(r,\rho,\eta)\in\boldsymbol{S}^{\mu,\nu}$ and $|\eta|\not=0$. Then 
\[ \op(a)(\eta):\Sch(\R\times X)\to\Sch(\R\times X) \]
extends to a continuous operator
\begin{equation} \label{cyleqqqF} 
\op(a)(\eta):H^{s,g}_\cone(X^\asymp)\to H^{s-\mu,g-\nu}_\cone(X^\asymp)
\end{equation}
for every $s,g\in\R$.
\end{Thm}
\proof
Let $u\in\Sch(\R\times X)$, and set $\|\cdot\|_{s,g}:=\|\cdot\|_{H^{s,g}_\cone(X^\asymp)}$, in particular, $\|\cdot\|_{0,0}=\|\cdot\|_{L^2(\R\times X)}$. By definition we have $\|u\|_{s,g}= \|\op(p^{s,g})(\eta^1)\|_{0,0}$.
Thus
\begin{multline*}
\|\op(a)(\eta)u\|_{s-\mu,g-\nu}=\|\op(p^{s-\mu,g-\nu})(\eta^1)\op(a)(\eta)u\|_{0,0} \\
=\|\op(p^{s-\mu,g-\nu})(\eta^1)\op(a)(\eta)\op(p^{s,g})^{-1}(\eta^1)\op(p^{s,g})(\eta^1)\|_{0,0}\leq c\|u\|_{s,g}
\end{multline*}
for $c:= \|\op(p^{s-\mu,g-\nu})(\eta^1)\op(a)(\eta)\op(p^{s,g})^{-1}(\eta^1)\|_{\Li(L^2(\R\times X))}$. It remains to prove that $c$ is a finite constant. This is completely straightforward when we replace $\op(p^{s,g})^{-1}(\eta^1)$ by $\op(p^{-s,-g})(\eta^1)$; in that case the remarks at the beginning of this section apply immediately, more precisely, we have
\begin{multline*} 
\op(p^{s-\mu,g-\nu})(\eta^1)\op(a)(\eta)\op(p^{-s,-g})(\eta^1)= \\ \op\big(p^{s-\mu,g-\nu} (\cdot,\eta^1)a(\cdot,\eta)p^{-s,-g}(\cdot,\eta^1)\big)
\end{multline*}
(where $\cdot$ stands for $r,\rho$) modulo a remainder of the form $\op(c)+R_N$ and $\op(c)$ is bounded in $L^2(\R\times X)$ for similar reasons as in Theorem \ref{cylCalVai} and the boundedness of $R_N$ in $L^2(\R\times X)$ is clear anyway. In general, $\op(p^{s,g})^{-1}(\eta^1)$ has the form $\op(p^{-s,-g})(\eta^1)+C_N(\eta^1)+R_N(\eta^1)$ for $C_N(\eta^1)=\op(c_N(\cdot,\eta^1))$ and a remainder $R_N(\eta^1)$ of smoothing behaviour while $c_N(\cdot,\eta)$ belongs to $\boldsymbol{S}^{-s-1,-g}$. Then, compared with the first step of the proof, we obtain extra terms, namely,
\begin{equation} \label{cyleeqqr1}
\op(p^{s-\mu,g-\nu})(\eta^1)\op(a)(\eta)\op(c_N)(\eta^1), 
\end{equation}
\begin{equation} \label{cyleeqqr2}
\op(p^{s-\mu,g-\nu})(\eta^1)\op(a)(\eta)R_N(\eta^1)
\end{equation}
which have to be bounded in $L^2(\R\times X)$. The arguments for \eqref{cyleeqqr1} are of the same structure as those at the beginning of the proof (the order of $c_N$ is even better than before), while for the composition \eqref{cyleeqqr2} we apply another remark from the preceding section, namely, that operators that are smoothing up to some degree (cf. Remark \ref{cylremmdeg}) when composed by other operators of the calculus give rise again to operators of sufficiently negative degree.
\qed

\end{document}